\documentclass[11pt,fancybox]{article}
\usepackage{amsmath,amssymb}
\usepackage{graphicx}
\newcommand{\dd}[1]{\langle #1\rangle}
\newcommand{\eproof}{\hspace*{\fill}{$\blacksquare$}}
\newtheorem{theorem}{Theorem}
\newtheorem{lemma}{Lemma}
\newtheorem{corollary}{Corollary}

\newtheorem{proposition}{Proposition}

\begin{document}
\title{
On canonical projections of the bipartite graph $D(k,q)$
\footnote{This work was supported by the National Natural Science Foundation of China (No. 61977056).}
\footnote{Email addresses: cz12024@163.com(C. Zhang), xycheng@yzu.edu.cn(X. Cheng), ystang@yzu.edu.cn(Y. Tang), mxu@szcu.edu.cn(M. Xu).}
}
\author{
Chen Zhang, Xiaoyan Cheng, Yuansheng Tang\footnote{Corresponding author.}\\
{\it\small School of Mathematical Sciences, Yangzhou University, Jiangsu, China}\\
and Ming Xu\\
{\it\small Suzhou City University,Jiangsu, China}
}

\date{}
\maketitle
\begin{abstract}

The bipartite graphs $D(k,q)$, constructed by Lazebnik et al. in 1995, are $q$-regular, edge-transitive, and of large girth.
In this paper, we primarily instigate the canonical projections of $D(k,q)$.
First, as a generalization of $D(k,q)$,
we define a bipartite graph $\Gamma(\Omega,q)$, where $\Omega$, referred to as a rooted binary tree, is a finite subset of binary sequences closed under the deletion of any suffix.
We then derive several significant properties regarding the symmetry, connectivity and girth of $\Gamma(\Omega,q)$. For arbitrary nonnegative integers $m$ and $n$, we fully characterize the connected components of $\Gamma(\Omega_{m,n},q)$, which is isomorphic to a canonical projection of $D(k,q)$ with $k\geq 2m+2$. Our results demonstrate that the connected components of $\{\Gamma(\Omega_{m,n},q)\}_{2\leq |m-n|\leq 3}$
form a novel family of graphs with large girth in the sense of Biggs.

{\it Index Terms}--Bipartite graph, canonical projection, rooted binary tree, automorphism, girth.

\end{abstract}

\section{Introduction}
All graphs we consider in this paper are assumed to be simple, i.e. undirected, without loops and multiple edges.
For a graph $G$, its vertex set and edge set are denoted by $V(G)$ and $E(G)$, respectively.
The {\it order} of $G$ is the number of vertices in $V(G)$.
The {\it size} of $G$ is the number of edges in $E(G)$.
The {\it degree} of a vertex of $G$ is the number of vertices adjacent to it.
A graph is {\it $r$-regular} if the degree of any vertex is equal to $r$.
A sequence $\omega=(v_1,v_1,\dots,v_n)$ of vertices in $V(G)$ is called a \emph{non-backtracking walk} (NB-walk) leading by $v_1$
if $v_{i+1}\sim v_i$ is an edge of $G$ for $1\leq i< n$ and $v_{i+1}\neq v_{i-1}$ for $1<i<n$.
The sequence $\omega$ is called a {\it closed non-backtracking walk} (CNB-walk) of length $n$ further if it is an NB-walk, $v_2\neq v_n\sim v_1\neq v_{n-1}$ and $n\geq 3$.
If $G$ is not a forest, then the {\it girth} of $G$, denoted by $g(G)$, is the minimum length of CNB-walks in $G$. Clearly, any CNB-walk of $G$ with length $g(G)$ is a {\it simple cycle}, which consists of distinct vertices.
A graph $G$ is {\it connected} if, for any two vertices $v,v'$ in $V(G)$, there is a walk from $v$ to $v'$.
For any graph $G$, let $c(G)$ denote the number of connected components of $G$.
A graph $G$ is \emph{edge-transitive} if, for any edges $v_1\sim v'_1$ and $v_2\sim v'_2$, there exists
an automorphism $\varphi\in \mathrm{Aut}(G)$ such that $\varphi(\{v_1,v_1'\})=\{v_2,v_2'\}$.

In literature, graphs with large girth and a high degree of symmetry
have been known to be hard to construct and have turned out to be
useful in different problems in extremal graph theory, finite geometry, coding theory, cryptography,
communication networks and quantum computations(c.f. \cite{K.P.P.P.F.2004}-\cite{Y.T.2011}).

Let $\mathcal{R}$ be an arbitrary commutative ring.
Some bipartite graphs defined by systems of equations over $\mathcal{R}$ were investigated in~\cite{L.W.2001}.
For $i\geq 2$, let $f_i: \mathcal{R}^{2i-2}\rightarrow \mathcal{R}$ be any given function.
For $k\geq 2$, Lazebnik et al. construct in~\cite{L.W.2001} a bipartite graph $B\Gamma_k=B\Gamma(\mathcal{R}; f_2,\ldots,f_k)$
whose vertex sets $L_k$ and $P_k$ are two copies of $\mathcal{R}^k$ such that two vertices $(l_1,l_2,\ldots,l_k)\in L_k$ and
$(p_1,p_2,\ldots,p_k)\in P_k$ are adjacent in $B\Gamma_k$ if and only if
\begin{align}
    l_i+p_i&=f_i(p_1,l_1,p_2,l_2,\ldots,p_{i-1},l_{i-1}), i=2,3,\ldots,k.
\end{align}
If $\mathcal{R}=\mathbb{F}_q$ is the finite field of
$q$ elements, $f_2(p_1,l_1)=p_1l_1$, $f_3(p_1,l_1,p_2,l_2)=p_1l_2$ and, for $4\leq i\leq k$,
\begin{align}
\label{dkq00}
f_i(p_1,l_1,p_2,l_2,\ldots,p_{i-1},l_{i-1})=\bigg\{
\begin{array}{ll}
l_1p_{i-2},&\text{if }i\equiv 0\text{ or }1\bmod 4,\\
p_1l_{i-2},&\text{if }i\equiv 2\text{ or }3\bmod 4,
\end{array}
\end{align}
then $B\Gamma_k$ is isomorphic to
the bipartite graph $D(k,q)$, introduced by Lazebnik et al. in \cite{L.U.1995},
which is $q$-regular, edge-transitive and of girth at least $k+4$.

For $2\leq l< k$, the {\it canonical projection} $\pi^k_l:L_k\cup P_k\to L_l\cup P_l$ is defined as the map
sending each vertex $v\in L_k\cup P_k$ to the vertex $v'\in L_l\cup P_l$ of the same type,
whose coordinates coincide with the first $l$ coordinates of $v$. Since $\pi^k_l$ is a {\it graph homomorphism} from $D(k,q)$ to $D(l,q)$ that maps adjacent vertices of $D(k,q)$ to adjacent vertices of $D(l,q)$, $D(l,q)$ is also referred to as a {\it canonical projection}
of $D(k,q)$.
Using these projections and edge transitivity, Lazebnik et al. fully characterized the connected components of $D(k,q)$ in \cite{L.U.W.1996} and \cite{L.V.2004}. Let $CD(k,q)$ be a connected component of $D(k,q)$. The graphs $CD(k,q)$ provide the best-known lower bounds for the maximum size of graphs with given order and girth, with the only exceptions being for girths 11 and 12, and for any prime power $q$ the graphs
$\{CD(k,q)\}_{k\geq 2}$ form a family of graphs with large girth in the sense of Biggs with constant $\gamma\geq\frac{4}{3}\log_q(q-1)$ (c.f. \cite{L.U.1995}, \cite{L.U.W.1995}).

In this paper, we focus on the canonical projections of $D(k,q)$ corresponding to the deletion of coordinates indexed by a small number of the largest even or odd integers.
We investigate their symmetry, connectivity and NB-walks in detail,
determine the exact number of connected components for each such projection,
and prove that a subset of these projections forms a new family of graphs that
matches the best-known lower bounds for the maximum size of graphs with given order and girth, just
as the original $D(k,q)$ graphs do.

The paper is organized as follows. A novel class of bipartite graphs $\Gamma(\Omega,q)$ over rooted binary trees $\Omega$ is constructed in Section~\ref{sec:02}. These graphs can be viewed as generalizations of the graphs $D(k,q)$ and
their canonical projections. The symmetry of $\Gamma(\Omega,q)$ is investigated in Section~\ref{sec:03}, where some automorphisms are constructed by extending the known automorphisms of $D(k,q)$.
The connectivity of $\Gamma(\Omega,q)$ is investigated in Section~\ref{sec:04},
where we construct several invariants for the connected components of $\Gamma(\Omega,q)$ and determine the exact number of connected components of $\Gamma(\Omega_{m,n},q)$, which is isomorphic to a canonical projection of the graph $D(k,q)$ with $k\geq 2m+2$.
For certain rooted binary trees $\Omega$, we investigate the NB-walks and CNB-walks of $\Gamma(\Omega,q)$ in Section~\ref{sec:05}, and then derive a lower bound for the girth of $\Gamma(\Omega,q)$. Finally, some concluding remarks are presented in Section~\ref{sec:06}.

\section{Bipartite Graphs over Rooted Binary Trees}
\label{sec:02}
In this section, we manage to generalize the definition of the graph $D(k,q)$ by indexing the coordinates of vertices with binary sequences.

Through this paper, we assume that $\Omega\subset \{0,1\}^*$ is a {\it rooted binary tree} (RB-tree) which is finite and contains all prefixes of any $\alpha\in\Omega$. Clearly, $\Omega$ always contains the empty sequence $\eta$, which is referred to as the {\it root} of $\Omega$. $\alpha\in\Omega$ is a {\it leaf} of $\Omega$ if $\{\alpha0,\alpha1\}\cap\Omega=\emptyset$.

Let $L(\Omega)$ and $R(\Omega)$ be two copies of
$\mathbb{F}_q^{|\Omega|+1}$. We denote vectors in
$L(\Omega)$ and $R(\Omega)$ by $[l]$ and $\dd{r}$, respectively, to distinguish their origins.
The coordinates of vectors in $L(\Omega)\cup R(\Omega)$ are indexed by
elements in $\Omega\cup\{*\}$, where $*$ is a symbol not in
$\{0,1\}^*$. For $[l]\in L(\Omega)$ and $\dd{r}\in R(\Omega)$, $l_*$ and $r_*$ are referred to as the {\it colors} of $[l]$ and $\dd{r}$, respectively.
Let $[0]_x$ denote the vector $[l]$ satisfying $l_*=x$ and $l_{\alpha}=0$ for all $\alpha\in\Omega$,
and let $\dd{0}_x$ denote the vector $\dd{r}$ satisfying $r_*=x$ and $r_{\alpha}=0$ for all $\alpha\in\Omega$.
Thus, $[0]_x$ and $\dd{0}_x$ are vectors that differ from the all-zero vector only in the coordinate indexed by $*$.

Let $\Gamma(\Omega,q)$ be the bipartite graph with vertex sets
$L(\Omega)$ and $R(\Omega)$ such that
$[l]\in L(\Omega)$ and $\dd{r}\in R(\Omega)$ are adjacent in
$\Gamma(\Omega,q)$ if
and only if
\begin{align}
\label{000}
    l_{\eta}+r_{\eta}=l_*r_*,
\end{align}
and, for any nonempty sequence $\alpha\in\Omega$,
\begin{align}\label{001}
&l_{\alpha}+r_{\alpha}=\bigg\{
\begin{array}{ll}
r_*l_{\beta},&\text{if }\alpha=\beta 0,\\
l_*r_{\beta},&\text{if }\alpha=\beta 1.
\end{array}
\end{align}
For convenience, we further define $*0=*1=\eta$, that is, the symbol $*$ removes the first bit of a subsequent sequence in $\{0,1\}^*$. Equation (\ref{000}) can thus be subsumed into equation (\ref{001}). It is clear that the bipartite graph
$\Gamma(\Omega,q)$ is $q$-regular, and the all-zero vectors
$[0]_0\in L(\Omega)$ and $\dd{0}_0\in R(\Omega)$ are adjacent in
$\Gamma(\Omega,q)$.

For $a\in\{0,1\}$ and $i\geq 1$,
let $\mu_{a,i}$ be the $(2i-1+a)$-th sequence
in the following set
\begin{align*}
U=\{0,1,01,10,010,101,0101,1010,\ldots\}.
\end{align*}
Clearly, $\mu_{a,i}$ is the $i$-th sequence in $U$ leading by $a$.
For convenience, we also define $\sigma_{a,0}=\eta$ and $\sigma_{a,-1}=*$ for $a\in\{0,1\}$.
For $m,n\geq 0$, let $\Omega_{m,n}$ be the RB-tree defined by
\begin{gather*}
\Omega_{m,n}=\{\mu_{0,i}:0\leq i\leq m\}\cup \{\mu_{1,j}:0\leq j\leq n\}.
\end{gather*}
Then, for $m,n\geq 0$, by interchanging the vertex sets we see easily
\begin{gather}
\Gamma(\Omega_{m,n},q)\cong\Gamma(\Omega_{n,m},q).
\label{yy0}
\end{gather}

Clearly, for $k\geq 2$ the graph $D(k,q)$ is isomorphic to a graph of form $\Gamma(\Omega_{m,n},q)$,
\begin{gather}
D(k,q)\cong\Gamma(\Omega_{\lceil (k-2)/2\rceil,\lfloor (k-2)/2\rfloor},q).
\label{yy1}
\end{gather}
Furthermore, for $0\leq n\leq m$, $\Gamma(\Omega_{m,n},q)$ is isomorphic to the canonical projection of the graph $D(2(m+1),q)$ corresponding to the deletion of coordinates with even indices that are greater than $2(n+1)$.

For the symmetry and girth of $D(k,q)$, the following proposition is from \cite{L.U.1995} and \cite{L.U.W.1995}.
\begin{proposition}\label{pp003}
Assume $q\geq 3$ is a prime power and $k\geq 2$. Then the graph $D(k,q)$ is edge-transitive and has girth at least $k+4$.
\end{proposition}

So far, the exact girth of $D(k,q)$ has been determined for many pairs $(k,q)$ (c.f. \cite{L.U.W.1995} and \cite{CCT.2014}-\cite{XCT23}).
For the connectivity of $D(k,q)$, the following proposition is from \cite{L.U.W.1996} and \cite{L.V.2004}.
\begin{proposition}\label{pp004}
Let $q\geq 3$ be a prime power, $k\geq 2$ be an integer, and $t=\lfloor(k+2)/4\rfloor$.
Then the connected components of $D(k,q)$ are isomorphic and amount to
\begin{enumerate}
\item $c(D(k,q))=q^{t-1}$ if $q \neq 4$,
\item $c(D(2,4))=c(D(3,4))=1$ and, $c(D(k,4))=4^t$ for $k\geq 4$.
\end{enumerate}
\end{proposition}

For the case $q=2$ and $k\geq 2$, it is shown in \cite{CCT.2014} that
the graph $D(k,2)$ consists of $2^{k-s-2}$ connected components, and each of them is a simple cycle of length $2^{s+3}$, where $s$ is the nonnegative integer satisfying
$2^{s+2}-3\leq k\leq 2^{s+3}-4$.

\section{Symmetry of $\Gamma(\Omega,q)$}\label{sec:03}
In this section, we show some automorphisms of $\Gamma(\Omega,q)$.
For any binary sequence $\alpha$, let $|\alpha|$ and $w(\alpha)$
denote its length and the number of its nonzero bits,
respectively.

The {\it multiplicative} automorphism of $\Gamma(\Omega,q)$ given in the following lemma, each coordinate is mapped into a multiple of itself,
is constructed by simply following \cite{L.U.W.1996} (or \cite{F.L.S.U.W.1995}).
\begin{lemma}
\label{lem:0001}
For $x,y\in \mathbb{F}_q^*$, let $\lambda_{x,y}$ denote the map
over $L(\Omega)\cup R(\Omega)$ such that, for $[l]\in L(\Omega)$
and $\dd{r}\in R(\Omega)$,
\begin{align*}
    (\lambda_{x,y}([l]))_*=xl_*,\ (\lambda_{x,y}(\dd{r}))_*=yr_*,
\end{align*}
and, for all $\alpha\in\Omega$,
\begin{align*}
    (\lambda_{x,y}([l]))_{\alpha}&=x^{w(\alpha)+1}y^{|\alpha|-w(\alpha)+1}l_{\alpha},\\
    (\lambda_{x,y}(\dd{r}))_{\alpha}&=x^{w(\alpha)+1}y^{|\alpha|-w(\alpha)+1}r_{\alpha}.
\end{align*}
Then, we have $\lambda_{x,y}\in\text{Aut}(\Gamma(\Omega,q))$.
\end{lemma}
{\bf Proof:} Clearly, $\lambda_{x,y}$ is a bijective map over
$L(\Omega)$ as well as over $R(\Omega)$ and,
$\lambda_{x,y}([l])\sim\lambda_{x,y}(\dd{r})$ if and only if $[l]\sim\dd{r}$. Hence, $\lambda_{x,y}$ belongs to $\text{Aut}(\Gamma(\Omega,q))$.
\eproof
\medskip

For the graph $D(k,q)$, a class of {\it additive} automorphisms was constructed in \cite{L.U.1995}, each of such automorphism adds an arbitrarily predesigned value to a specific coordinate, and then $D(k,q)$ was proved to be edge-transitive and of girth $\geq k+4$. Now we manage to generalize these additive automorphisms for some of the graphs $\Gamma(\Omega,q)$.
Let $\mathcal{S}(\Omega)=\cup_{\alpha\in\Omega}S_{\Omega}(\alpha)$, where $S_{\Omega}(\alpha)$ denotes the set of sequences $\beta\in\{0,1\}^*$ with
$\{\alpha01\beta,\alpha10\beta\}\cap\Omega\neq\emptyset$.

\begin{lemma}
\label{lem:0002} If $\alpha\in\Omega$ satisfies $S_{\Omega}(\alpha)\subset\Omega$,
then for any $x\in \mathbb{F}_q$ there is an automorphism
$\theta_{x,\alpha}\in\text{Aut}(\Gamma(\Omega,q))$ such that, for any
$[l]\in L(\Omega)$ and $\dd{r}\in R(\Omega)$,
\begin{align}
    (\theta_{x,\alpha}([l]))_{\alpha}=l_{\alpha}+x,\
    (\theta_{x,\alpha}(\dd{r}))_{\alpha}=r_{\alpha}-x \label{400},
\end{align}
and
\begin{align}
\label{401}
    (\theta_{x,\alpha}([l]))_{\gamma}=l_{\gamma},\
    (\theta_{x,\alpha}(\dd{r}))_{\gamma}=r_{\gamma},
\end{align}
for all $\gamma\in\Omega_{\alpha}=(\{*\}\cup\Omega)
\setminus\big(\{\alpha,\alpha0,\alpha1\}\cup
\{\alpha01\beta,\alpha10\beta:\beta\in\{0,1\}^*\}\big)$.
\end{lemma}
{\bf Proof:} Assume $x\in \mathbb{F}_q$. Let $\theta_{x,\alpha}$ be the map defined
over $L(\Omega)\cup R(\Omega)$ satisfying (\ref{400}),
\begin{gather*}
    (\theta_{x,\alpha}([l]))_{\alpha0}=l_{\alpha0}\text{ and } (\theta_{x,\alpha}(\dd{r}))_{\alpha0}=r_{\alpha0}+xr_*\text{ if }\alpha0\in \Omega,\\
    (\theta_{x,\alpha}([l]))_{\alpha1}=l_{\alpha1}-xl_*\text{ and } (\theta_{x,\alpha}(\dd{r}))_{\alpha1}=r_{\alpha1}\text{ if }\alpha1\in \Omega,\\
    (\theta_{x,\alpha}([l]))_{\alpha01\beta}=l_{\alpha01\beta}+xl_{\beta}\text{ and } (\theta_{x,\alpha}(\dd{r}))_{\alpha01\beta}=r_{\alpha01\beta}+xr_{\beta}\text{ if }\alpha01\beta\in \Omega,\\
    (\theta_{x,\alpha}([l]))_{\alpha10\beta}=l_{\alpha10\beta}-xl_{\beta}\text{ and } (\theta_{x,\alpha}(\dd{r}))_{\alpha10\beta}=r_{\alpha10\beta}-xr_{\beta}\text{ if }\alpha10\beta\in \Omega,
\end{gather*}
and (\ref{401}) for $\gamma\in\Omega_{\alpha}$.
Clearly, $\theta_{x,\alpha}$ is well-defined over $L(\Omega)\cup
R(\Omega)$ and a bijective map over $L(\Omega)$ as well as over
$R(\Omega)$. One can also check easily that
$\theta_{x,\alpha}([l])\sim\theta_{x,\alpha}(\dd{r})$ if and only if $[l]\sim\dd{r}$. Hence, the map $\theta_{x,\alpha}$ is indeed
the desired automorphism of $\Gamma(\Omega,q)$.
\eproof
\medskip

For any binary sequence $\alpha\in\{0,1\}^*$ and $i\geq 0$, let $\alpha^i$ denote the
following sequence
\begin{align*}
    \alpha^i=\bigg\{
    \begin{array}{ll}
    \underbrace{\alpha\alpha\cdots \alpha}_i\ ,&\text{ if }i>0,\\
    \eta,&\text{ if }i=0.
    \end{array}
\end{align*}

\begin{theorem}
\label{thm:-002}
Assume that $\mathcal{S}(\Omega)\subset\Omega$.
Then, for any edge $[l]\sim\dd{r}$ of $\Gamma(\Omega,q)$,
\begin{enumerate}
    \item there is an automorphism $\theta_0\in\text{Aut}(\Gamma(\Omega,q))$ such that $\theta_0([l])=[0]_{l_*}$ and
\begin{align}
    &(\theta_0(\dd{r}))_*=r_*, \label{703}\\
    &(\theta_0(\dd{r}))_{1^i}=l_*^{i+1}r_*, \text{ if }1^i\in\Omega, i\geq 0,\label{704}\\
    &(\theta_0(\dd{r}))_{\alpha}=0,\text{ for all }\alpha\in\Omega\setminus\{1^i:i\geq 0\},\label{705}
\end{align}
    \item there is an automorphism $\theta_1\in\text{Aut}(\Gamma(\Omega,q))$ such that $\theta_1(\dd{r})=\dd{0}_{r_*}$ and
\begin{align}
    &(\theta_1([l]))_*=l_*, \\
    &(\theta_1([l]))_{0^i}=r_*^{i+1}l_*, \text{ if }0^i\in\Omega, i\geq 0,\\
    &(\theta_1([l]))_{\alpha}=0,\text{ for all }\alpha\in\Omega\setminus\{0^i:i\geq 0\}.
\end{align}
\end{enumerate}
Furthermore,
the components of $\Gamma(\Omega,q)$ are isomorphic.
\end{theorem}
{\bf Proof:}
Since, for any $\alpha\in\Omega$, the
automorphism $\theta_{x,\alpha}$ given in Lemma~\ref{lem:0002}
fixes those coordinates of $[l]\in L(\Omega)$ and $\dd{r}\in
R(\Omega)$ whose indices are in $\{*\}\cup\{\beta\in \Omega:
|\beta|\leq |\alpha|,\beta\neq \alpha\}$, we see that there is an
automorphism $\theta_0\in\text{Aut}(\Gamma(\Omega,q))$ satisfying
$\theta_0([l])=[0]_{l_*}$ and (\ref{703}). Then, (\ref{704}) and
(\ref{705}) follow from
$\theta_0([l])\sim\theta_0(\dd{r})$ and
(\ref{001}). Hence, we have proved the first conclusion. The second conclusion can be proved similarly.

Let $C$ be an arbitrary component of $\Gamma(\Omega,q)$.
Clearly, $C$ must contain a vertex $[l]$ with color $l_*=0$. Then, according to the first conclusion, there exists an
automorphism $\theta_0\in\text{Aut}(\Gamma(\Omega,q))$ such that $\theta_0([l])=[0]_{0}$.
Hence, $C$ is isomorphic to the component which contains the vertex $[0]_0$.
\eproof
\medskip

Clearly, if $\mathcal{S}(\Omega)\subset\Omega$, according to Lemma~\ref{lem:0002} and Theorem~\ref{thm:-002} we see that any vertex $[l]\in L(\Omega)$ with $l_*\neq 0$ may be mapped to $[0]_1$ by some automorphism of $\Gamma(\Omega,q)$.

Since for any nonnegative integers $m,n$ we have
\begin{gather}
\mathcal{S}(\Omega_{m,n})\subset\Omega_{m,n}
\text{ if and only if }|m-n|\leq 3,
\end{gather}
the following corollary follows from Theorem~\ref{thm:-002} directly.
\begin{corollary}
\label{coro002}
Let $m,n$ be nonnegative inters with $|m-n|\leq 3$.
Then,
the components of $\Gamma(\Omega_{m,n},q)$ are isomorphic. In particular, for any edge $[l]\sim\dd{r}$ of $\Gamma(\Omega_{m,n},q)$,
\begin{enumerate}
    \item there is an automorphism $\theta_0\in\text{Aut}(\Gamma(\Omega_{m,n},q))$ such that $\theta_0([l])=[0]_{l_*}$ and
\begin{gather*}
    (\theta_0(\dd{r}))_*=r_*, \
    (\theta_0(\dd{r}))_{\eta}=l_*r_*, \\
    (\theta_0(\dd{r}))_1=l^2_*r_*, \text{ if }n\geq 1,\\
    (\theta_0(\dd{r}))_{\alpha}=0,\text{ for all }\alpha\in\Omega_{m,n}\setminus\{\eta,1\},
\end{gather*}
    \item there is an automorphism $\theta_1\in\text{Aut}(\Gamma(\Omega_{m,n},q))$ such that $\theta_1(\dd{r})=\dd{0}_{r_*}$ and
\begin{gather*}
    (\theta_1([l]))_*=l_*, \ (\theta_1([l]))_{\eta}=l_*r_*, \\
    (\theta_1([l]))_0=l_*r_*^2, \text{ if }m\geq 1,\\
    (\theta_1([l]))_{\alpha}=0,\text{ for all }\alpha\in\Omega_{m,n}\setminus\{\eta,0\}.
\end{gather*}
\end{enumerate}
\end{corollary}

\section{Connectivity of $\Gamma(\Omega,q)$}\label{sec:04}
In this section, we mainly discuss the connectivity of the graph $\Gamma(\Omega,q)$.
At first we show some invariants for the components of $\Gamma(\Omega,q)$ by following the methods shown in \cite{L.U.W.1996} and \cite{L.V.2004}. Then, we reveal some significant relations between the components of $\Gamma(\Omega,q)$ and those of their canonical projections. The exact number of components of $\Gamma(\Omega_{m,n},q)$ is determined in the end of this section.
We always assume that $[l]\sim\dd{r}$ is an edge of $\Gamma(\Omega,q)$ without specification
in this section.

The invariants contained in the following lemma are essentially translations of those shown in  \cite{L.U.W.1996}, a new proof for them is given here.
\begin{lemma}
\label{lem:d000}
For any $\alpha,\beta\in\{0,1\}^*\cup\{*\}$, we have
\begin{gather}
\label{y001} r_{\beta}l_{\alpha1}-r_{\alpha}l_{\beta1}
=r_{\alpha}r_{\beta1}
-r_{\beta}r_{\alpha1}, \text{ if }\{\alpha 1, \beta 1\}\subset\Omega,
\\
\label{y002} l_{\beta}r_{\alpha0}-l_{\alpha}r_{\beta0}
=l_{\alpha}l_{\beta0}
-l_{\beta}l_{\alpha0},\text{ if }\{\alpha 0, \beta 0\}\subset\Omega,\\
\label{y100} r_{\beta}l_{\alpha01}-r_{\beta10}l_{\alpha}
=(r_{\alpha0}r_{\beta1}-r_{\beta}r_{\alpha01})+(l_{\alpha}l_{\beta10}
-l_{\beta1}l_{\alpha0}),\text{ if }\{\alpha01,\beta10\}\subset\Omega.
\end{gather}
Furthermore, if $k>0$ and $\{(01)^k,(10)^k\}\subset\Omega$, then
\begin{gather}
\psi_k([l])=l_{(10)^k}-l_{(01)^k}+\sum_{j=0}^{k-1}\big(l_{(01)^{k-1-j}}l_{(10)^j}-
l_{*(01)^j}l_{*(10)^{k-j}}
\big),\label{y101}\\
\psi_k(\dd{r})=r_{(01)^k}-r_{(10)^k}+
\sum_{i=0}^{k-1}\big(r_{(10)^{k-1-i}}r_{(01)^i}-r_{*(10)^i}r_{*(01)^{k-i}}\big)
\label{y101'}
\end{gather}
invariants for the components of $\Gamma(\Omega,q)$. Indeed, for each edge $[l]\sim\dd{r}$ of $\Gamma(\Omega,q)$ we have $\psi_k([l])=\psi_k(\dd{r})$.
\end{lemma}
{\bf Proof:}
If $\{\alpha 1, \beta 1\}\subset\Omega$, from
(\ref{001}) we see
\begin{align*}
r_{\beta}l_{\alpha1}-r_{\alpha}l_{\beta1}
=r_{\beta}(l_*r_{\alpha}-r_{\alpha1})-r_{\alpha}(l_*r_{\beta}-r_{\beta1})
=r_{\alpha}r_{\beta1}-r_{\beta}r_{\alpha1},
\end{align*}
i.e., we have (\ref{y001}).
Similarly, we have (\ref{y002}). Therefore, from (\ref{y001})
and (\ref{y002}) we see
\begin{align*}
&r_{\beta}l_{\alpha01}-r_{\beta10}l_{\alpha}\\
=&r_{\beta}l_{\alpha01}-r_{\alpha0}l_{\beta1}
+l_{\beta1}r_{\alpha0}-l_{\alpha}r_{\beta10}\\
=&(r_{\alpha0}r_{\beta1}-r_{\beta}r_{\alpha01})+(l_{\alpha}l_{\beta10}
-l_{\beta1}l_{\alpha0}).
\end{align*}
Hence, we have (\ref{y100}).

Furthermore, if $\{(01)^k,(10)^k\}\subset\Omega$, according to (\ref{001}) and (\ref{y100}) we see
\begin{align*}
&\big(l_{(10)^k}+r_{(10)^k}\big)-\big(l_{(01)^k}+r_{(01)^k}\big)
=r_*l_{(10)^{k-1}1}-l_*r_{(01)^{k-1}0}\\
=&r_*l_{*(01)^k}-r_{*(10)^k}l_*=
\sum_{i=0}^{k-1}\big(r_{*(10)^i}l_{*(01)^{k-i}}-r_{*(10)^{i+1}}l_{*(01)^{k-i-1}}\big)\\
=&\sum_{i=0}^{k-1}\Big(\big(r_{*(01)^{k-i-1}0}r_{*(10)^i1}-r_{*(10)^i}r_{*(01)^{k-i}}\big)
+\big(l_{*(01)^{k-i-1}}l_{*(10)^{i+1}}
-l_{*(10)^i1}l_{*(01)^{k-i-1}0}\big)\Big)\\
=&\sum_{i=0}^{k-1}\Big(\big(r_{(10)^{k-i-1}}r_{(01)^i}-r_{*(10)^i}r_{*(01)^{k-i}}\big)
+\big(l_{*(01)^{k-i-1}}l_{*(10)^{i+1}}
-l_{(01)^i}l_{(10)^{k-i-1}}\big)\Big).
\end{align*}
Then, according to (\ref{y101}) and (\ref{y101'}) we see $\psi_k([l])=\psi_k(\dd{r})$ for each edge $[l]\sim\dd{r}$ and thus
$\psi_k([l])$ and $\psi_k(\dd{r})$ are invariants for the components of $\Gamma(\Omega,q)$.
\eproof
\medskip

The invariants shown in the following lemma can be found in \cite{L.V.2004}, though therein only $\psi_0([l])=\psi_0(\dd{r})=0$ was used to characterize the component containing the vertex $[0]_0$.
\begin{lemma}\label{lem:d001}
Assume $q=4$ and $\{0,1\}\subset\Omega$. Then, $\psi_0([l])=l_*l_{\eta}+l_1+l_0^2$ and $\psi_0(\dd{r})=r_*r_{\eta}+r_0+r_1^2$ are invariants for the components of $\Gamma(\Omega,4)$.
Indeed, for any edge $[l]\sim\dd{r}$ of $\Gamma(\Omega,4)$ we have
\begin{gather}
\label{200}\psi_0([l])=\psi_0(\dd{r})^2,\\
\label{201}\psi_0(\dd{r})=\psi_0([l])^2.
\end{gather}
\end{lemma}
{\bf Proof:} Since for any $x\in\mathbb{F}_4$ we have $x^4=x$, according to (\ref{001}) we see
\begin{align*}
\psi_0([l])+\psi_0(\dd{r})^2&=l_*l_{\eta}+l_1+l_0^2+(r_*r_{\eta}+r_0+r_1^2)^2\\
&=l_*l_{\eta}+(l_1+r_1)+(r_*r_{\eta}+r_0+l_0)^2\\
&=l_*(l_{\eta}+r_{\eta})+(r_*r_{\eta}+r_*l_{\eta})^2\\
&=l_*^2r_*+(r_*^2l_*)^2=0,
\end{align*}
and thus we have (\ref{200}) and (\ref{201}). Subsequently, $\psi_0([l])$ and $\psi_0(\dd{r})$ are invariants for the components of $\Gamma(\Omega,4)$.
\eproof\medskip

Now we suppose that $\Omega'\subset\Omega$ is also an RB-tree.
Let $\pi_{\Omega/\Omega'}$ denote the {\it canonical projection} from $L(\Omega)\cup R(\Omega)$ to $L(\Omega')\cup R(\Omega')$
defined by deleting the entries indexed by sequences in $\Omega\setminus\Omega'$.
By convention, we assume $\pi_{\Omega/\Omega'}(L(\Omega))=L(\Omega')$ and $\pi_{\Omega/\Omega'}(R(\Omega))=R(\Omega')$. For any bipartite graph $G$ with
$V(G)\subseteq L(\Omega)\cup R(\Omega)$, let $G'$ be the bipartite graph with
$V(G')\subseteq L(\Omega')\cup R(\Omega')$ such that $u', v'\in V(G')$
are adjacent in $G'$ if and only if there exist vertices $u\in \pi^{-1}_{\Omega/\Omega'}(u')$ and $v\in\pi^{-1}_{\Omega/\Omega'}(v')$ that are adjacent in $G$.
The graph $G'$ is referred to the {\it canonical projection} of $G$ under $\pi_{\Omega/\Omega'}$, and written as $G'=\pi_{\Omega/\Omega'}(G)$.
For a component $C$ of $\Gamma(\Omega,q)$ and a vertex $u\in V(\pi_{\Omega/\Omega'}(C))$, let $\nu(u,C)$ denote the set of vertices $v\in V(C)$ with $\pi_{\Omega/\Omega'}(v)=u$.

The following lemma is a generalization of Propositions 3.3 and 3.4 in \cite{L.U.W.1996}.
\begin{lemma}\label{gg01}
For any component $C$ of $\Gamma(\Omega,q)$, $C'=\pi_{\Omega/\Omega'}(C)$ is a component of $\Gamma(\Omega',q)$ and
$\pi_{\Omega/\Omega'}$ is a $t$-to-1 graph homomorphism from $C$ to $C'$ for some $t$ with $1\leq t\leq q^{|\Omega|-|\Omega'|}$. Furthermore, for each tree $T'$ in $C'$ there are just pairwise vertex-disjoint trees $T_1,T_2,\ldots,T_t$ in $C$ such that $\pi_{\Omega/\Omega'}(T_i)=T'$, $i=1,2,\ldots,t$, and
$\pi_{\Omega/\Omega'}$ is an isomorphism from $C$ to $C'$ if $t=1$.
\end{lemma}
{\bf Proof}:
Since for any two adjacent vertices $v, v'$ of $\Gamma(\Omega,q)$, their projections $\pi_{\Omega/\Omega'}(v), \pi_{\Omega/\Omega'}(v')$ are adjacent in $\Gamma(\Omega',q)$,
it is clearly that $C'$ is a connected subgraph of $\Gamma(\Omega',q)$.
If $u\in V(C')$ and $u'\in V(\Gamma(\Omega',q))$ are adjacent in $\Gamma(\Omega',q)$, it is clear that, for any $v\in \nu(u,C)$,
there is a unique vertex $v'\in V(\Gamma(\Omega,q))$ with $\pi_{\Omega/\Omega'}(v')=u'$
such that $v'$ and $v$ are adjacent in $\Gamma(\Omega,q)$. Hence, we have $v'\in V(C)$ and thus $u'\in V(C')$.
This implies that $C'$ is a component of $\Gamma(\Omega',q)$.

If $u,u'\in V(C')$ are adjacent in $C'$ and $\nu(u,C)=\{v_1,v_2,\ldots,v_{t_u}\}$, $t_u=|\nu(u,C)|$, then there are vertices $v_1',v_2',\ldots,v_{t_u}'\in\nu(u',C)$
such that $v_i'$ and $v_i$ are adjacent in $\Gamma(\Omega,q)$ for $i=1,2,\ldots,t_u$.
Clearly, the vertices $v_1',v_2',\ldots,v_{t_u}'$ are distinct. Hence, we have $t_{u'}\geq t_u$ and thus $t=t_u$ is independent of $u$.
Clearly, we see $1\leq t\leq q^{|\Omega|-|\Omega'|}$ and that $\pi_{\Omega/\Omega'}$ is a $t$-to-1 graph homomorphism from $C$ to $C'$.

Furthermore, we see easily that, for each tree $T'$ in $C'$, there are just pairwise vertex-disjoint trees $T_1,T_2,\ldots,T_t$ in $C$ such that $\pi_{\Omega/\Omega'}(T_i)=T'$, $i=1,2,\ldots,t$, and
$\pi_{\Omega/\Omega'}$ is an isomorphism from $C$ to $C'$ if $t=1$.
\eproof
\medskip

Let $\alpha$ be a leaf of $\Omega$ and set $\Omega'=\Omega\setminus\{\alpha\}$. Clearly, $\Omega'$ is an RB-tree.
For any component $C$ of $\Gamma(\Omega,q)$ and $u\in V(\pi_{\Omega/\Omega'}(C))$, let $\rho(u,C)$ denote the set
of $v_{\alpha}$ for $v\in\nu(u,C)$, i.e.
\begin{align*}
\rho(u,C)=\bigg\{
\begin{array}{ll}
\{l_{\alpha}:[l]\in \nu(u,C)\},&\text{if }u\in L(\Omega'),\\
\{r_{\alpha}:\dd{r}\in \nu(u,C)\},&\text{if }u\in R(\Omega').
\end{array}
\end{align*}
For any component $C'$ of $\Gamma(\Omega',q)$, let $\mathcal{L}(C')$ denote the set of
components $C$ of $\Gamma(\Omega,q)$ with $\pi_{\Omega/\Omega'}(C)=C'$.
The following theorem is a refinement of Lemma~\ref{gg01} for the case $\Omega'=\Omega\setminus\{\alpha\}$.

\begin{theorem}\label{gg02}
Suppose that $\alpha$ is a leaf of $\Omega$ and $C'$ is a component of $\Gamma(\Omega',q)$, where $\Omega'=\Omega\setminus\{\alpha\}$.
Then, $s=|\mathbb{L}(C')|$ divides $q$ and there exist maps
$f: V(C')\to \mathbb{F}_q$, $g: \mathbb{L}(C')\to \mathbb{F}_q$ and an additive subgroup $G$ of $\mathbb{F}_q$ with order
$t=q/s$ such that
\begin{enumerate}
  \item For any $C\in\mathbb{L}(C')$ and $v\in V(C')$,
\begin{gather}
\rho(v,C)=G+f(v)+s(v)g(C)\label{cc03}
\end{gather}
is a coset of $G$,
where $s(v)=1$ if $v\in L(\Omega')$ and $s(v)=-1$ if $v\in R(\Omega')$;
  \item $\{g(C):C\in \mathbb{L}(C')\}$ is a representive set of cosets of $G$ in $\mathbb{F}_q$, namely, $\{g(C)+G\}_{C\in\mathbb{L}(C')}$ is a partition of $\mathbb{F}_q$.
\end{enumerate}
In particular, the components in $\mathbb{L}(C')$ are isomorphic.
Furthermore, if $\mathcal{S}(\Omega')\subset\Omega'$,
then the subgroup $G$ is independent of components of $\Gamma(\Omega',q)$ and
$t=|G|=p^{mn}$ for some integer $m$ with $0\leq m\leq \log_{p^{n}}q$, where $p$ is the characteristic of $\mathbb{F}_q$ and
$n$ is the smallest positive integer such that $q-1$ divides $(p^{n}-1)\gcd(q-1,w(\alpha)+1,|\alpha|-w(\alpha)+1)$.
\end{theorem}
{\bf Proof}:
Clearly, for any $x\in\mathbb{F}_q$, the automorphism $\theta_{x,\alpha}$ given in Lemma~\ref{lem:0002} defines a bijective map in each
$\mathcal{L}(C')$ and
\begin{align}
\label{t300}
\rho(u,\theta_{x,\alpha}(C))=
\rho(u,C)+s(u)x\text{ for any }u\in V(C')\text{ and }C\in\mathcal{L}(C').
\end{align}

Since for any vertex $u\in V(C')$ the set $\{\rho(u,C):C\in\mathcal{L}(C')\}$ forms a partition of $\mathbb{F}_q$,
there exists a unique component in $\mathcal{L}(C')$, denoted by $C_u$, such that $0\in \rho(u,C_u)$.
Hence, for any $u\in V(C')$ and $y\in \rho(u,C_u)$, according to (\ref{t300}) we have
$y\in\rho\big(u,\theta_{s(u)y,\alpha}(C_u)\big)\cap\rho(u,C_u)\neq\emptyset$
and thus we see
\begin{gather}
\label{t303}
\theta_{s(u)y,\alpha}(C_u)=C_u,\\
\label{t301}
\rho(u,C_u)=\rho(u,C_u)+s(u)y.
\end{gather}
From (\ref{t301}) and the arbitrariness of $y$ we see that $\rho(u,C_u)$ is an additive subgroup of $\mathbb{F}_q$.

Furthermore, from (\ref{t300}) and (\ref{t303}) we see that
$\rho(v,C_u)=\rho(v,C_u)+y$ is valid for any $v\in V(C')$ and $y\in \rho(u,C_u)$.
Hence, $\rho(v,C_u)$ is a union of several cosets of $\rho(u,C_u)$ in $\mathbb{F}_q$.
Since according to Lemma~\ref{gg01} the cardinality of $\rho(v,C_u)$ is independent of $v\in V(C')$,
we see that the set $\rho(v,C_u)$ is just a coset of $\rho(u,C_u)$ in $\mathbb{F}_q$.

Since for any $C\in\mathcal{L}(C')$ we have $\{\theta_{x,\alpha}(C):x\in\mathbb{F}_q\}=\mathcal{L}(C')$, for any vertices $u,v\in V(C')$ there is an element $x\in \mathbb{F}_q$ such that $\theta_{x,\alpha}(C_u)=C_{v}$ and then
according to (\ref{t300}) we see $\rho(v,C_{v})=\rho(v,\theta_{x,\alpha}(C_u))=\rho(v,C_u)+s(v)x$ is also a coset of $\rho(u,c_u)$. Hence,
the group $G=\rho(u,C_u)$ is independent of $u\in V(C')$. Let $t$ be the order of $G$. Since $G$ is a subgroup of $\mathbb{F}_q$, $t$ is a factor of $q$ and $s=|\mathcal{L}(C')|=q/t$ divides $q$.

Let $u$ be a given vertex in $V(C')$. For any $C\in\mathcal{L}(C')$ and adjacent vertices $v,v'$ of $C'$, let $a,b$ be elements in $\mathbb{F}_q$ such that $\rho(v,C)=\rho(v,C_u)+a$ and $\rho(v',C)=-\rho(v,C)+b$, where $b$ depends only on $v,v'$ and is independent of $C$.
Then,
\begin{gather*}
\rho(v',C)=-(\rho(v,C_u)+a)+b=(-\rho(v,C_u)+b)-a
=\rho(v',C_u)-a.
\end{gather*}
Hence, there is an element $g(C)\in\mathbb{F}_q$ such that
$\rho(v,C)=\rho(v,C_u)+s(v)g(C)$ for any $v\in V(C')$.
Therefore, for the map $f: V(C')\to \mathbb{F}_q$ with $\rho(v,C_u)=G+f(v)$, $v\in V(C')$,
we have (\ref{cc03}).
Since $\{\rho(v,C):C\in \mathbb{L}(C')\}$ is a partition of $\mathbb{F}_q$, $\{g(C):C\in \mathbb{L}(C')\}$ forms a representive set of cosets of $G$ in $\mathbb{F}_q$.

For $v\in V(C')$ and $c\in \mathbb{F}_q$, let $u=u(v,c)$ denote the vertex in $\pi_{\Omega/\Omega'}^{-1}(v)$ with $u_{\alpha}=c$.
For components $C,\tilde{C}$ in $\mathcal{L}(C')$, let $\kappa:V(C)\to V(\tilde{C})$ be the map defined by $\kappa(u(v,c))=u\big(v,c+s(v)(g(\tilde{C})-g(C))\big)$ for any $v\in V(C')$ and $c\in\rho(v,C)$. Then, $\kappa$ must be an isomorphism from $C$ to $\tilde{C}$.
Hence, the components in $\mathcal{L}(C')$ are isomorphic.

Now we assume $\mathcal{S}(\Omega')\subset\Omega'$. Let $\beta$ be an arbitrary sequence in $\Omega$ and $x$ an arbitrary element in $\mathbb{F}_q$.
From $\mathcal{S}(\Omega')\subset\Omega'=\Omega\setminus\{\alpha\}$, we see $\{\beta01\alpha,\beta10\alpha\}\cap\Omega=\emptyset$.
Then,
for the automorphism $\theta_{x,\beta}\in\text{Aut}(\Gamma(\Omega,q))$ and any vertex $v$ of $\Gamma(\Omega,q)$,
according to the proof of Lemma~\ref{lem:0002} one can confirm easily that $\pi_{\Omega/\Omega'}(\theta_{x,\beta}(v))$ and
$(\theta_{x,\beta}(v))_{\alpha}-v_{\alpha}$ are determined only by the projection $\pi_{\Omega/\Omega'}(v)$ of $v$.
Hence, for any component $C$ of $\Gamma(\Omega,q)$ and vertices $v,v'\in V(C)$ with $\pi_{\Omega/\Omega'}(v)=\pi_{\Omega/\Omega'}(v')$, we have $(\theta_{x,\beta}(v))_{\alpha}-(\theta_{x,\beta}(v'))_{\alpha}=v_{\alpha}-v'_{\alpha}$
and thus the subgroup $G$ given for the component $\pi_{\Omega/\Omega'}(C)$ is the same one for the component $\pi_{\Omega/\Omega'}(\theta_{x,\beta}(C))$.
Since $\alpha$ is a leaf of $\Omega$, we have $S(\alpha)=\emptyset$ and thus $\mathcal{S}(\Omega)=\mathcal{S}(\Omega')\subset\Omega'\subset\Omega$.
Then, according to the proof of Theorem~\ref{thm:-002}, we see that any component of $\Gamma(\Omega,q)$ is an image of $C$ under a few automorphisms $\theta_{x,\beta}\in\text{Aut}(\Gamma(\Omega,q))$,
and thus it can be concluded further that the subgroup $G$ is independent of the components of $\Gamma(\Omega',q)$.

Now we assume further that
the component $C'$ contains $u_0$, one of the all-zero vertices of $\Gamma(\Omega',q)$.
Let $\xi$ be a primitive element of $\mathbb{F}_q$. Since for any integers $i,j$ the multiplicative automorphism $\lambda_{\xi^i,\,\xi^j}$ fixes the vertex $u_0\in V(C')$ and the component $C_{u_0}\in\mathbb{L}(C')$, we see that
\begin{align}\label{cc11}
    \xi^{i(w(\alpha)+1)+j(|\alpha|-w(\alpha)+1)}\rho(u_0,C_{u_0})=\rho(u_0,C_{u_0})\text{ for any }i,j.
\end{align}
Since $G=\rho(u_0,C_{u_0})$ is an additive subgroup of $\mathbb{F}_q$,
we see that (\ref{cc11}) is equivalent to that $xG=G$ is valid for any nonzero element $x$ in
\begin{align*}
    \Delta=\bigg\{\sum_{0\leq i<(q-1)/k}a_i\xi^{ki}: a_i\in\mathbb{F}_p\bigg\},
\end{align*}
where $k=\gcd(q-1,w(\alpha)+1,|\alpha|-w(\alpha)+1)$ and $p$ is the characteristic of $\mathbb{F}_q$.
Clearly, the set $\Delta$ is the smallest field which contains $\xi^k$ and thus
$\Delta=\mathbb{F}_{p^n}$, where $n>0$ is the smallest integer such that $(q-1)|(p^{n}-1)k$.
Hence, $G$ can be seen as a subspace of $\mathbb{F}_q$ over $\mathbb{F}_{p^n}$ and thus there is an integer $m$ with $0\leq m\leq \log_{p^n}q$
such that $t=|G|=p^{mn}$.
\eproof
\medskip

For example, for $a\in\{0,1\}$ and $\alpha=\mu_{a,2i+1}$, we have $\gcd(w(\alpha)+1,|\alpha|-w(\alpha)+1)=1$ and thus the power $p^n$ in Theorem~\ref{gg02} is equal to $q$, this implies $\Delta=\mathbb{F}_q$ and that $G$ is either $\mathbb{F}_q$ or $\{0\}$.
 For any prime power $q\geq 3$, the number of components of $\Gamma(\Omega_{m,n},q)$ is determined in the following theorem.
\begin{theorem}
\label{theo:401}
Let $q\geq 3$ be a prime power.
For $m\geq n\geq 0$, $c(\Gamma(\Omega_{m,n},q))$ is equal to
\begin{gather}
\label{h001}c_{n,q}=\bigg\{
\begin{array}{ll}
q^{\lfloor n/2\rfloor},&\text{if }q\neq 4\text{ or }n=0,\\
4^{\lfloor n/2\rfloor+1},&\text{if }q= 4\text{ and }n\geq 1.
\end{array}
\end{gather}
\end{theorem}
{\bf Proof}: The proof proceeds by induction on the difference $m-n$.

If $0\leq n\leq m\leq n+1$, according to (\ref{yy1}) and Proposition \ref{pp004} we see the number of components of the graph $\Gamma(\Omega_{m,n},q)\cong D(m+n+2,q)$ is equal to $q^{\lfloor (m+n)/4\rfloor}=q^{\lfloor n/2 \rfloor}=c_{n,q}$ if $q\neq 4$ or $n=0$, and to $4^{\lfloor (m+n)/4\rfloor+1}=4^{\lfloor n/2\rfloor+1}=c_{n,4}$ if $q=4$ and $n\geq 1$.
Therefore, according to (\ref{yy0}) we have $c(\Gamma(\Omega_{m,n},q))=c_{n,q}$ when $0\leq m-n\leq 1$.

Assume $c(\Gamma(\Omega_{m,n},q))=c_{n,q}$ is valid when $0\leq m-n\leq t$. For given $n\geq 0$, let $\Omega=\Omega_{n+t+1,n+1}$, $\Omega'=\Omega_{n+t+1,n}$ and $\Omega''=\Omega_{n+t,n}$.
The inductive step is then treated in three separate cases.

Case 1: $q\neq 4$ and $2|n$, or $q=4$ and $n\geq 2$ is even. From (\ref{h001}) we see
\begin{gather*}
c(\Gamma(\Omega,q))=c_{n+1,q}=c_{n,q}=c(\Gamma(\Omega'',q)).
\end{gather*}
According to $\Gamma(\Omega'',q)=\pi_{\Omega'/\Omega''}(\Gamma(\Omega',q))$, $\Gamma(\Omega',q)=\pi_{\Omega/\Omega'}(\Gamma(\Omega,q))$ and Lemma~\ref{gg01},
we see
\begin{gather*}
c(\Gamma(\Omega'',q))\leq c(\Gamma(\Omega',q))\leq c(\Gamma(\Omega,q)).
\end{gather*}
Hence, we have $c(\Gamma(\Omega_{n+t+1,n},q))=c(\Gamma(\Omega',q))=c_{n,q}$.

Case 2: $q=4$ and $n=0$. Clearly, we have $\{0,1\}\subset\Omega$, $\Omega\setminus\Omega'=\{1\}$ and
\begin{gather*}
c(\Gamma(\Omega,4))=c(\Gamma(\Omega_{t+1,1},4))=c_{1,4}=4.
\end{gather*}
According to Lemma~\ref{lem:d001}, $\psi_0([l])$ is invariant in any component $C$ of $\Gamma(\Omega,4)$. Since the value of $\psi_0([l])-l_1$ is determined by the entries of $[l']=\pi_{\Omega/\Omega'}([l])$,
we see $|\rho([l'],C)|=1$ and thus the group $G$ given in Theorem~\ref{gg02} for $C'=\pi_{\Omega/\Omega'}(C)$ is $\{0\}$. Hence, we have $|\mathcal{L}(C')|=4$ and thus we see
\begin{align*}
c(\Gamma(\Omega_{t+1,0},4))=c(\Gamma(\Omega',4))=4^{-1}c(\Gamma(\Omega,4))=1=c_{0,4}.
\end{align*}

Case 3: $n=2k-1$ is odd. Clearly, we have $\{(01)^k,(10)^k\}\subset\Omega$,
$\Omega\setminus\Omega'=\{(10)^k\}$ and
\begin{gather*}
c(\Gamma(\Omega,q))=c(\Gamma(\Omega_{n+t+1,n+1},q))=c_{n+1,q}=qc_{n,q}.
\end{gather*}
According to Lemma~\ref{lem:d000}, $\psi_k([l])$ is invariant in any component $C$ of $\Gamma(\Omega,q)$.
Since the value of $\psi_k([l])-l_{(10)^k}$ is determined by the entries of $[l']=\pi_{\Omega/\Omega'}([l])$,
we see $|\rho([l'],C)|=1$ and thus the group $G$ given in Theorem~\ref{gg02} for $C'=\pi_{\Omega/\Omega'}(C)$ is $\{0\}$.
Hence, we have $|\mathcal{L}(C')|=q$ and thus we see
\begin{align*}
c(\Gamma(\Omega_{n+t+1,n},q))=c(\Gamma(\Omega',q))= q^{-1}c(\Gamma(\Omega,q))=c_{n,q}.
\end{align*}

Therefore, we have proved $c(\Gamma(\Omega_{m,n},q))=c_{n,q}$ when $m-n=t+1$.
\eproof
\medskip

Clearly, according the proof of Theorem~\ref{theo:401} and the argument in \cite{L.V.2004}, we see easily that the connected components of $\Gamma(\Omega_{m,n},q)$
can be fully characterized by the invariants shown in Lemmas~\ref{lem:d000} and \ref{lem:d001}.
In particular, for $m\geq 1$ each component of $\Gamma(\Omega_{m,1},4)$ is isomorphic to $\Gamma(\Omega_{m,0},4)$, and for $m\geq 2k\geq 2$ any component $C$ of $\Gamma(\Omega_{m,2k},q)$ is isomorphic to $\pi_{\Omega_{m,2k}/\Omega_{m,2k-1}}(C)$.
\section{Girth of $\Gamma(\Omega,q)$}\label{sec:05}

To give a lower bound for the girth of certain graphs $\Gamma(\Omega,q)$,
 in this section we mainly investigate the NB-walks in $\Gamma(\Omega,q)$ whose leading vertices are of form $[0]_x$ or $\dd{0}_x$.
\subsection{NB-walks in $\Gamma(\Omega,q)$}\label{sec:051}
In this subsection, we consider to express the vertices on an NB-walk by
their colors if the beginning vertex is of form $[0]_x$. According to symmetry, similar argument is also true for the case that the NB-walks of $\Gamma(\Omega,q)$
start at vertices of form $\dd{0}_x$.

Let
$[l^{(1)}]\sim\dd{r^{(2)}}\sim[l^{(3)}]\sim\dd{r^{(4)}}\sim\cdots$
be a given NB-walk
in $\Gamma(\Omega,q)$. Assume $(x_1,x_2,\ldots)$ is the sequence of colors of the vertices in this NB-walk, namely, $x_i=l^{(i)}_*$ if $i$ is odd, and $x_i=r^{(i)}_*$ if $i$ is even.
Clearly, according to the definition of NB-walks we have $x_i\neq x_{i+2}$ for all $i\geq 1$.
For $j\geq 1$, $s\geq 0$ and
$\alpha\in\{0,1\}^*\cup\{*\}$, from $[l^{(2j-1)}]\sim\dd{r^{(2j)}}\sim [l^{(2j+1)}]\sim\dd{r^{(2j+2)}}$ and (\ref{001}) we see
\begin{gather}
\label{801}
    l^{(2j+1)}_{\alpha10^s}-l^{(2j-1)}_{\alpha10^s}=x_{2j}^s(x_{2j+1}-x_{2j-1})r^{(2j)}_{\alpha},
    \text{ if }\alpha10^s\in\Omega,\\
\label{802}
    r^{(2j+2)}_{\alpha01^s}-r^{(2j)}_{\alpha01^s}=x_{2j+1}^s(x_{2j+2}-x_{2j})l^{(2j+1)}_{\alpha},
    \text{ if }\alpha01^s\in\Omega,
\end{gather}
where $x^0=1$ for all $x\in\mathbb{F}_q$.

For positive integer $i$ and sequence $\boldsymbol{s}=(s_1,s_2,\ldots)$ of
nonnegative integers, let $\sigma_i(\boldsymbol{s})$ denote the sequence defined iteratively by
\begin{gather}
\label{t00'}
\sigma_{i}(\boldsymbol{s})=
\bigg\{
\begin{array}{ll}
\sigma_{i-1}(\boldsymbol{s})10^{s_{i}},&\text{if $i$ is even},\\
\sigma_{i-1}(\boldsymbol{s})01^{s_{i}},&\text{if $i$ is odd},
\end{array}
\end{gather}
where $\sigma_1(\boldsymbol{s})=*$ if $s_1=0$ and $\sigma_1(\boldsymbol{s})=1^{s_1-1}$
if $s_1>0$.

Clearly, for any sequence $\alpha\in\{0,1\}^*$,
there uniquely exist some nonnegative integers $s_1,s_2,\ldots,s_{2n}$ such that
\begin{align}
\alpha=1^{s_1}0^{s_2+1}1^{s_3+1}0^{s_4+1}\cdots1^{s_{2n-1}+1}0^{s_{2n}}
=\sigma_{2n}(s_1,s_2,\ldots,s_{2n}),
\end{align}
which can also be expressed as
\begin{align}
\alpha=\bigg\{\begin{array}{ll}
\sigma_{2n-1}(s_1,s_2,\ldots,s_{2n-2},s_{2n-1}+1),&\text{if }s_{2n}=0,\\
\sigma_{2n+1}(s_1,s_2,\ldots,s_{2n-1},s_{2n}-1,0),&\text{if }s_{2n}>0.
\end{array}
\end{align}

Hereafter, let $\boldsymbol{s}=(s_1,s_2,\ldots)$ be a given sequence of nonnegative integers.
\begin{lemma}
If $[l^{(1)}]=[0]_{x_1}$,
then for any positive integer $i$ we have
\begin{gather}
\label{a01}
l^{(2i-1)}_{\sigma_{2j}(\boldsymbol{s})}=0, \text{ if }j\geq i\text{ and }\sigma_{2j}(\boldsymbol{s})\in\Omega,\\
\label{a02}
r^{(2i)}_{\sigma_{2j+1}(\boldsymbol{s})}=0, \text{ if }j\geq i\text{ and }\sigma_{2j+1}(\boldsymbol{s})\in\Omega.
\end{gather}
\end{lemma}
{\bf Proof}:
From $[l^{(1)}]=[0]_{x_1}$, we see that (\ref{a01}) is valid for $i=1$.
Furthermore, for $j\geq 1$ from $[l^{(1)}]\sim\dd{r^{(2)}}$ and (\ref{801}) to (\ref{t00'}) we have
\begin{align*}
r^{(2)}_{\sigma_{2j+1}(\boldsymbol{s})}
=&r^{(2)}_{\sigma_{2j}(\boldsymbol{s})01^{s_{2j+1}}}+l^{(1)}_{\sigma_{2j}(\boldsymbol{s})01^{s_{2j+1}}}
=x_1^{s_{2j+1}}r^{(2)}_{\sigma_{2j}(\boldsymbol{s})0}
\\
=&x_1^{s_{2j+1}}\Big(r^{(2)}_{\sigma_{2j}(\boldsymbol{s})0}
+l^{(1)}_{\sigma_{2j}(\boldsymbol{s})0}\Big)
=x_1^{s_{2j+1}}x_2l^{(1)}_{\sigma_{2j}(\boldsymbol{s})}=0.
\end{align*}
Hence, (\ref{a02}) is valid for $i=1$.

Now we assume that (\ref{a01}) and (\ref{a02}) are valid for $1\leq i\leq t$.
For $j\geq t$, from (\ref{801}) and (\ref{t00'}) we see
\begin{align}
l^{(2t+1)}_{\sigma_{2j+2}(\boldsymbol{s})}
=l^{(2t+1)}_{\sigma_{2j+2}(\boldsymbol{s})}-l^{(2t-1)}_{\sigma_{2j+2}(\boldsymbol{s})}
=(x_{2t+1}-x_{2t-1})x_{2t}^{s_{2j+2}}r^{(2t)}_{\sigma_{2j+1}(\boldsymbol{s})}
=0.\label{a04}
\end{align}
From (\ref{802}) and (\ref{t00'}) we see
\begin{align*}
r^{(2t+2)}_{\sigma_{2j+3}(\boldsymbol{s})}
=r^{(2t+2)}_{\sigma_{2j+3}(\boldsymbol{s})}-r^{(2t)}_{\sigma_{2j+3}(\boldsymbol{s})}
=(x_{2t+2}-x_{2t})x_{2t+1}^{s_{2j+3}}l^{(2t+1)}_{\sigma_{2j+2}(\boldsymbol{s})}
=0.
\end{align*}
Hence, (\ref{a01}) and (\ref{a02}) are valid for $i=t+1$.

Then, by induction we have (\ref{a01}) and (\ref{a02}) for all $i\geq 1$.\eproof\medskip

For positive integers $a_1,a_2,\ldots$ with
$a_{2j-1}\leq a_{2j}<a_{2j+1}$, $j=1,2,\ldots,$
let
\begin{align}
\label{a60}
\Phi_{\boldsymbol{s}}(a_1,a_2,\ldots,a_{n})
=&\bigg\{
\begin{array}{ll}
\Phi_{\boldsymbol{s}}(a_1,a_2,\ldots,a_{n-1})
(x_{2a_{n}+1}-x_{2a_{n}-1})x_{2a_{n}}^{s_{n}},&\text{if $n$ is even},\\
\Phi_{\boldsymbol{s}}(a_1,a_2,\ldots,a_{n-1})
(x_{2a_{n}}-x_{2a_{n}-2})x_{2a_{n}-1}^{s_{n}},&\text{if $n$ is odd},
\end{array}
\end{align}
where $\Phi_{\boldsymbol{s}}(a_1)=(x_{2a_1}-x_{2a_1-2})x_{2a_1-1}^{s_1}$ and $x_0=0$.

\begin{theorem}
If $[l^{(1)}]=[0]_{x_1}$, then, for positive integers $i$, $j$ with $1\leq j\leq i$,
\begin{align}
\label{a06}
r^{(2i)}_{\sigma_{2j-1}(\boldsymbol{s})}
&=\sum_{1\leq a_1\leq a_2<a_3\leq a_4<\cdots<a_{2j-1}\leq i}
\Phi_{\boldsymbol{s}}(a_1,a_2,\ldots,a_{2j-1}),
\text{ if }\sigma_{2j-1}(\boldsymbol{s})\in\Omega,\\
\label{a07}
l^{(2i+1)}_{\sigma_{2j}(\boldsymbol{s})}
&=\sum_{1\leq a_1\leq a_2<a_3\leq a_4<\cdots<a_{2j-1}\leq a_{2j}\leq i}
\Phi_{\boldsymbol{s}}(a_1,a_2,\ldots,a_{2j}),\text{ if }\sigma_{2j}(\boldsymbol{s})\in\Omega.
\end{align}
In particular, for $i\geq 1$ we have
\begin{align}
\label{a26}
r^{(2i)}_{\sigma_{2i-1}(\boldsymbol{s})}&=x_2x_1^{s_1}\prod_{a=2}^{2i-1}(x_{a+1}-x_{a-1})x_a^{s_{a}},
\text{ if }\sigma_{2i-1}(\boldsymbol{s})\in\Omega,\\
\label{a27}
l^{(2i+1)}_{\sigma_{2i}(\boldsymbol{s})}&=x_2x_1^{s_1}\prod_{a=2}^{2i}(x_{a+1}-x_{a-1})x_a^{s_{a}},
\text{ if }\sigma_{2i}(\boldsymbol{s})\in\Omega.
\end{align}
\end{theorem}
{\bf Proof}:
For $1\leq j\leq i$, from (\ref{801}), (\ref{t00'}) and(\ref{a01}), we have
\begin{align}
\label{a10}
l^{(2i+1)}_{\sigma_{2j}(\boldsymbol{s})}
=\sum_{a=j}^i\big(l^{(2a+1)}_{\sigma_{2j}(\boldsymbol{s})}-l^{(2a-1)}_{\sigma_{2j}(\boldsymbol{s})}
\big)
=\sum_{a=j}^i(x_{2a+1}-x_{2a-1})x_{2a}^{s_{2j}}r^{(2a)}_{\sigma_{2j-1}(\boldsymbol{s})}.
\end{align}
For $2\leq j\leq i$, from (\ref{802}), (\ref{t00'}) and (\ref{a02}), we have
\begin{align}
\label{a11}
r^{(2i)}_{\sigma_{2j-1}(\boldsymbol{s})}
=\sum_{a=j}^i\big(r^{(2a)}_{\sigma_{2j-1}(\boldsymbol{s})}
-r^{(2a-2)}_{\sigma_{2j-1}(\boldsymbol{s})}\big)
=\sum_{a=j}^i(x_{2a}-x_{2a-2})x_{2a-1}^{s_{2j-1}}l^{(2a-1)}_{\sigma_{2j-2}(\boldsymbol{s})}.
\end{align}
If $s_1>0$, we have $\sigma_1(\boldsymbol{s})=1^{s_1-1}$, and thus
from $[l^{(1)}]=[0]_{x_1}\sim\dd{r^{(2)}}$ and (\ref{001}) we see
\begin{align*}
r^{(2)}_{\sigma_1(\boldsymbol{s})}=x_2x_1^{s_1}=(x_2-x_0)x_1^{s_1}=\Phi_{\boldsymbol{s}}(1),
\end{align*}
and then, for $i>1$, from $\sigma_1(\boldsymbol{s})=*01^{s_1-1}$ and (\ref{802}) we have
\begin{align*}
r^{(2i)}_{\sigma_1(\boldsymbol{s})}
=&r^{(2)}_{\sigma_1(\boldsymbol{s})}+\sum_{a=2}^i
\left(r^{(2a)}_{\sigma_1(\boldsymbol{s})}-r^{(2a-2)}_{\sigma_1(\boldsymbol{s})}\right)\\
=&(x_2-x_0)x_1^{s_1}+\sum_{a=2}^i(x_{2a}-x_{2a-2})x_{2a-1}^{s_{1}-1}l^{(2a-1)}_{*}\\
=&\sum_{a=1}^i(x_{2a}-x_{2a-2})x_{2a-1}^{s_{1}}=\sum_{a=1}^i\Phi_{\boldsymbol{s}}(a).
\end{align*}
If $s_1=0$, then $\sigma_1(\boldsymbol{s})=*$ and from $x_0=0$ we see
\begin{align*}
r^{(2i)}_{\sigma_1(\boldsymbol{s})}=x_{2i}=\sum_{a=1}^i(x_{2a}-x_{2a-2})x_{2a-1}^{s_{1}}
=\sum_{a=1}^i\Phi_{\boldsymbol{s}}(a).
\end{align*}
Hence, (\ref{a06}) is valid for $j=1$.
Therefore, from (\ref{a10}), (\ref{a11}) and the definition of $\Phi_{\boldsymbol{s}}(\cdot)$,
according to induction on $j$ we see (\ref{a06}) and (\ref{a07}).

Clearly, according to (\ref{a06}) and (\ref{a07}) we have
\begin{gather*}
r^{(2i)}_{\sigma_{2i-1}(\boldsymbol{s})}=\Phi_{\boldsymbol{s}}(1,1,2,2,\ldots,i-1,i-1,i),\\
l^{(2i+1)}_{\sigma_{2i}(\boldsymbol{s})}=\Phi_{\boldsymbol{s}}(1,1,2,2,\ldots,i,i),
\end{gather*}
and thus from (\ref{a60}) we see (\ref{a26}) and (\ref{a27}).
\eproof\medskip

\subsection{CNB-walks and Girth of $\Gamma(\Omega,q)$}\label{sec:052}

In this subsection, we will show some conditions for the existence of
some CNB-walks which contain a vertex of form $[0]_x$ or $\dd{0}_x$.
Furthermore, some lower bounds for the girth of $\Gamma(\Omega,q)$ are deduced
according to these conditions.

Assume that
\begin{gather*}
[l^{(1)}]=[0]_{x_1}\sim\dd{r^{(2)}}\sim\cdots\sim[l^{(2k-1)}]\sim\dd{r^{(2k)}},
\end{gather*}
denoted by $\omega$, is an arbitrary
CNB-walk in $\Gamma(\Omega,q)$ of length $2k$.
Clearly, the colors of the vertices in $\omega$ satisfy
\begin{gather}
x_i\neq x_{i+2},\text{ for all }1\leq i\leq 2k,\label{a30}
\end{gather}
where $x_{2k+1}=x_1$ and $x_{2k+2}=x_2$.

\begin{lemma}
\label{lem:200}
If $m>0$ and $\mu_{0,m}\in\Omega$, then $\Gamma(\Omega,q)$ has no CNB-walk of length $2(m+1)$ containing a vertex of form $[0]_x$.
\end{lemma}
{\bf Proof}:
At first, we assume that $m=2i-1$ is odd. Then, $\mu_{0,m}=(01)^{i-1}0$. Assume that the CNB-walk $\omega$ is of length $2k=4i=2(m+1)$.
Since
\begin{gather*}
[l^{(1)}]\sim\dd{r^{(2)}}\sim[l^{(3)}]\sim\cdots\sim\dd{r^{(2i)}}\sim[l^{(2i+1)}],\\
[l^{(1)}]\sim\dd{r^{(4i)}}\sim[l^{(4i-1)}]\sim\cdots\sim\dd{r^{(2i+2)}}\sim[l^{(2i+1)}]
\end{gather*}
are two NB-walks of length $2i+1$, for $s\in\{0,1\}$ from
$(01)^{i-1}0^s=\sigma_{2i}(0,\ldots,0,s)\in\Omega$ and (\ref{a27}) we have
\begin{align*}
    l^{(2i+1)}_{(01)^{i-1}0^s}&=x_{2i}^sx_2(x_3-x_1)\prod_{a=3}^{2i}(x_{a+1}-x_{a-1})\\
                       &=x_{2i+2}^sx_{4i}(x_{4i-1}-x_1)\prod_{a=3}^{2i}(x_{4i-a+1}-x_{4i-a+3}).
\end{align*}
If $l^{(2i+1)}_{(01)^{i-1}}=0$, then we have $x_2=0=x_{4i}=x_{2k}$, contradicts (\ref{a30}).
If $l^{(2i+1)}_{(01)^{i-1}}\neq0$, then we have
\begin{align*}
x_{2i}=l^{(2i+1)}_{(01)^{i-1}0}\left/ l^{(2i+1)}_{(01)^{i-1}}=x_{2i+2}\right.,
\end{align*}
still contradicts (\ref{a30}).

Now we assume that $m=2j$ is even. Then, $\mu_{0,m}=(01)^{j}$. Assume that the CNB-walk $\omega$ is of length $2k=4j+2=2(m+1)$.
Since
\begin{gather*}
[l^{(1)}]\sim\dd{r^{(2)}}\sim[l^{(3)}]\sim\cdots\sim[l^{(2j+1)}]\sim\dd{r^{(2j+2)}},\\
[l^{(1)}]\sim\dd{r^{(4j+2)}}\sim[l^{(4j+1)}]\sim\cdots\sim[l^{(2j+3)}]\sim\dd{r^{(2j+2)}}
\end{gather*}
are two NB-walks of length $2j+2$, for $s\in\{0,1\}$ from
$0(10)^{j-1}1^s=\sigma_{2j+1}(0,\ldots,0,s)\in\Omega$ and (\ref{a26}) we have
\begin{align*}
    r^{(2j+2)}_{0(10)^{j-1}1^s}&=x_{2j+1}^sx_2(x_3-x_1)\prod_{a=3}^{2j+1}(x_{a+1}-x_{a-1})\\
                       &=x_{2j+3}^sx_{4j+2}(x_{4j+1}-x_1)\prod_{a=3}^{2j+1}(x_{4j-a+3}-x_{4j-a+5}).
\end{align*}
If $r^{(2j+2)}_{0(10)^{j-1}}=0$, then we have $x_2=0=x_{4j+2}=x_{2k}$, contradicts (\ref{a30}).
If $r^{(2j+2)}_{0(10)^{j-1}}\neq0$, then we have
\begin{align*}
x_{2j+1}=r^{(2j+2)}_{0(10)^{j-1}1}\left/ r^{(2j+2)}_{0(10)^{j-1}}=x_{2j+3}\right.,
\end{align*}
still contradicts (\ref{a30}).
\eproof
\medskip

At the end of this subsection, we give a lower bound for girth of some graphs $\Gamma(\Omega_{m,n},q)$.
\begin{theorem}\label{theo:502}
For nonnegative integers $m,n$ with $2\leq |m-n|\leq 3$, the graph $\Gamma(\Omega_{m,n},q)$ has girth at least $2(\max\{m,n\}+2)$.
\end{theorem}
{\bf Proof}:
Since $\Gamma(\Omega_{m,n},q)\cong\Gamma(\Omega_{n,m},q)$, we assume without loss of generality that $m>n$.
Let $\omega$ be an arbitrary CNB-walk of length $2k$ in $\Gamma(\Omega_{m,n},q)$.
According to Corollary~\ref{coro002}, we can assume further that $\omega$ contains a vertex
of form $[0]_x$.
Therefore, from Lemma~\ref{lem:200} and $\mu_{0,i}\in\Omega_{m,n}$ for $1\leq i\leq m$, we see that $k\not\in\{2,3,\ldots,m+1\}$ and thus we have $k\geq m+2$.
Hence, $\Gamma(\Omega_{m,n},q)$ has girth at least $2(m+2)$.
\eproof
\medskip

According to Propositions~\ref{pp003}, \ref{pp004} and Theorems~\ref{theo:401}, \ref{theo:502}, for any $l\geq 0$ we see that the connected components of $\Gamma(\Omega_{2l+3,2l},q)$ and $D(4l+5,q)\cong\Gamma(\Omega_{2l+2,2l+1},q)$ have identical order, size and lower bound on their girth. In particular, the connected components of $\{\Gamma(\Omega_{m,n},q)\}_{2\leq|m-n|\leq3}$ form a novel family of graphs with large girth in the sense of Biggs.

\section{Conclusion}\label{sec:06}

To investigate the canonical projections of $D(k,q)$, this paper introduces a new class of bipartite graphs, denoted as $\Gamma(\Omega,q)$, where $\Omega\in\{0,1\}^*$ is a rooted-binary tree containing all prefixes of each $\alpha\in\Omega$. We analyze the symmetry, connectivity and girth of these graphs. For certain instances of $\Gamma(\Omega,q)$, we construct several automorphisms and invariants for their connected components.
For any NB-walk $\omega$ in $\Gamma(\Omega,q)$ leading by $[0]_x$ or $\dd{0}_x$, we derive explicit formulas to express the coordinates of vertices in $\omega$ in terms of their colors.
For any graph $\Gamma(\Omega_{m,n},q)$, which is isomorphic to a canonical projection of $D(k,q)$ with $k\geq 2m+2$,
we determine the exact number of its connected components. For nonnegative integers $m,n$ satisfying $2\leq |m-n|\leq 3$, we establish a nontrivial lower bound for the girth of $\Gamma(\Omega_{m,n},q)$. Our results demonstrate that the graphs $\Gamma(\Omega_{m,n},q)$ under the condition $2\leq |m-n|\leq 3$
perform comparably to $D(k,q)$ in several respects, in particular, the connected components of $\{\Gamma(\Omega_{m,n},q)\}_{2\leq|m-n|\leq3}$ form a novel family of graphs with large girth in the sense of Biggs.

\end{document}